\documentclass[10pt]{article}
\usepackage{times}
\usepackage{amsthm}
\usepackage{amssymb}
\usepackage{latexsym}
\usepackage{amscd}
\usepackage{epsfig}
\renewcommand{\labelenumi}{(\theenumi)}

\newtheorem{thm}{Theorem}[section]
\newtheorem{lem}[thm]{Lemma}
\newtheorem{que}[thm]{Question}

\newcommand{\R}{\ensuremath{\mathbf{R}}}

\newcommand{\RRR}{\ensuremath{\mathbf{R}^3}}

\newcommand{\Z}{\ensuremath{\mathbf{Z}}}
\newcommand{\bd}{\ensuremath{\partial}}
\newcommand{\irr}{ir\-re\-duc\-i\-ble}

\newcommand{\birr}{\bd-\irr}

\newcommand{\Inte}{Int \, }
\newcommand{\nul}{\emptyset}
\newcommand{\p}{^{\prime}}
\newcommand{\pp}{^{\prime\prime}}
\newcommand{\bp}{^{\, \prime}}

\begin{document}

\title{\normalsize \textbf{UNCOUNTABLY MANY ARCS IN} $S^3$  
\textbf{WHOSE COMPLEMENTS HAVE \\ NON-ISOMORPHIC,  
INDECOMPOSABLE FUNDAMENTAL GROUPS}}
\author{\small ROBERT MYERS \\ \small \textit{Department of Mathematics} 
\\ \small \textit{Oklahoma State University} \\  
\small \textit{Stillwater, OK 74078 USA} \\ 
\small \textit{Email: myersr@math.okstate.edu}
\vspace{.7in}}
\date{}

\maketitle
\begin{abstract} 
An uncountable collection of arcs in $S^3$ is 
constructed, each member of which is wild precisely at its endpoints, 
such that the fundamental groups of their complements are 
non-trivial, pairwise non-isomorphic, and indecomposable 
with respect to free products.
The fundamental group of the complement of a certain    
Fox-Artin arc is also shown to be indecomposable. 

\vspace{\baselineskip}


\end{abstract}

\section*{\normalsize 1. Introduction}
\stepcounter{section}
At the 1996 Workshop in Geometric Topology F. D. Ancel \cite{An} 
posed the following questions:

\begin{que} Let $A$ be the Fox-Artin arc 
in $S^3$ which is pictured in Figure 1. 
Is $\pi_1(S^3-A)$ indecomposable with respect to free products? 
\end{que}

\begin{figure}[!h]
\begin{center}
\epsfig{file=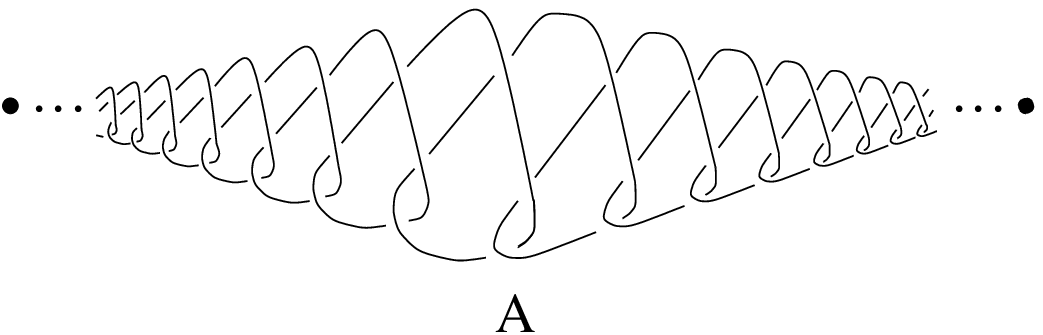, width=4in}
\end{center}
\caption{The Fox-Artin arc $A$}
\end{figure}

\begin{que} Are there infinitely (uncountably?) many wild arcs 
$A_i$ in $S^3$ such that $\pi_1(S^3-A_i)$ and $\pi_1(S^3-A_j)$ are 
non-isomorphic for $i \neq j$? \end{que}

Fox and Artin \cite{FA} proved that $\pi_1(S^3-A)$ is non-trivial. 
($A$ is actually the mirror image of their Example 1.1.) 
At the workshop Ancel remarked that an incorrect proof that it is 
indecomposable had been published by Ros\l aniec 
\cite{Ro}. He also noted that an affirmative answer to 
Question 1.1 would give an affirmative answer to the countable case of 
Question 1.2 by concatenating finitely many copies of $A$; the resulting 
groups are free products of copies of $\pi_1(S^3-A)$ and so would be  
non-isomorphic \cite[Vol. II, p. 27]{Ku}. 
These examples would have a finite but unbounded number of wild points. 

In this paper we answer these two questions in the affirmative. 
In particular, regarding Question 1.2 we construct an uncountable 
family of arcs $A_i$ such that the fundamental groups $\pi_1(S^3-A_i)$ 
are non-isomorphic for distinct indices and also are indecomposable 
and non-trivial. Moreover each arc is wild precisely 
at its endpoints. 

We remark that if the fundamental group of the complement 
of an arc in $S^3$ is non-trivial, then it is not finitely 
generated \cite[Corollary 2.6]{GHM}. 

Ancel also posed the following question, to which one can of course 
add the question of indecomposability. As of this writing these  
questions remain open, but it seems likely 
that affirmative answers could be obtained by the methods of this 
paper. 

\begin{que} Let $B$ be the wild arc in the solid torus $V$ pictured 
in Figure 2. Suppose $k_i:V \rightarrow S^3$ is a knotted embedding 
such that $\pi_1(S^3-k_i(V))$ is not isomorphic to $\pi_1(S^3-k_j(V))$ 
for $i \neq j$. Is $\pi_1(S^3-k_i(B))$ not isomorphic to   
$\pi_1(S^3-k_j(V))$ for $i \neq j$? \end{que}

\begin{figure}[h]
\begin{center}
\epsfig{file=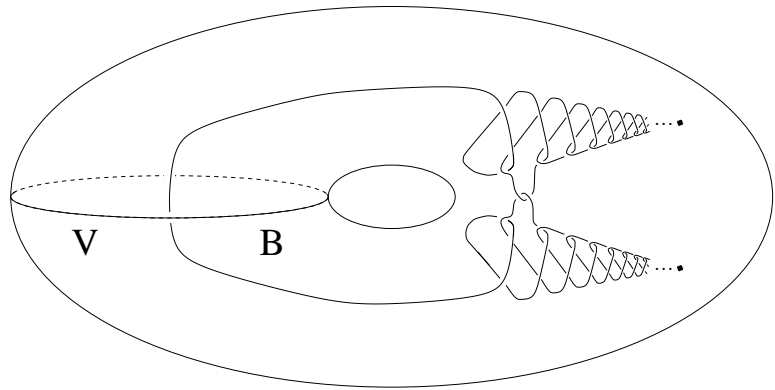, width=4in}
\end{center}
\caption{Tie $V$ in a knot to get an Ancel arc.}
\end{figure}

The paper is organized as follows. 
In section 2 we give a criterion for 
the fundamental group of a non-compact 3-manifold to be indecomposable 
and non-trivial. 
In section 3 we prove that the exterior of the Fox-Artin arc satisfies 
this criterion. 
In section 4 we prove a lemma about 
embeddings of torus knot groups in torus knot groups. In section 5 
we construct the uncountable family of arcs mentioned above and 
verify its properties. 

The author thanks Bill Banks for drawing the Fox-Artin arc which is 
used in Figures 1, 2, and 3. 

\section*{\normalsize 2. A Criterion for Indecomposability}
\stepcounter{section}

Recall that a group $G$ is \textit{decomposable} if it is a free product 
$K*L$, where $K$ and $L$ are non-trivial. $G$ is \textit{indecomposable} 
if it is not decomposable. 

\begin{lem} Let $\{H_k\}_{k \geq 0}$ be a sequence of 
non-trivial, non-infinite-cyclic, indecomposable 
subgroups of $G$ such that $H_k 
\subseteq H_{k+1}$ for all $k \geq 0$ and $G= \cup^{\infty}_{k=0} H_k$. 
Then $G$ is indecomposable. \end{lem} 

\begin{proof} Suppose $G=K*L$, where $K$ and $L$ are non-trivial. 
Then no non-trivial element of $K$ is conjugate to an element of $L$. 
This can be seen as follows. Let $N$ be the normal closure of $K$ 
in $G$. Let $p:G \rightarrow G/N$ be the natural projection. Then 
there is an isomorphism $q:G/N \rightarrow L$ such that the restriction 
of $q \circ p$ to $L$ is the identity of $L$ \cite[pp. 101--102]{Ma}. 
But $q \circ p$ sends any conjugate of an element of $K$ to the trivial 
element of $L$. 

By the Kurosh subgroup theorem \cite{Ku,Ma} any subgroup 
of $G$ is 
a free product of a free group and conjugates of subgroups of $K$ and 
of $L$. Since $H_0$ in indecomposable and non-infinite-cyclic we may 
thus assume that it is conjugate to a subgroup of $K$. Similarly $H_1$ 
must be conjugate to a subgroup of $K$ or of $L$. The latter cannot 
happen since then some non-trivial element of $K$ would be conjugate  
to an element of $L$. Continuing in this fashion we get that each $H_k$ 
is conjugate to a subgroup of $K$. This implies that $G$ cannot be the 
union of the $H_k$ since the non-trivial elements of $L$ are excluded. 
$\square$
\end{proof}

We now consider fundamental groups of non-compact 3-manifolds. For 
basic definitions in 3-manifold topology we refer to \cite{He} and  
\cite{Ja}. A 3-manifold $M$ is \textit{\birr} if $\bd M$ is 
incompressible in $M$. Let $S$ and $S\p$ be compact surfaces such that 
$S$ is properly embedded in $M$ and $S\p$ either is properly embedded 
in $M$ or lies in $\bd M$. Then $S$ and $S\p$ are \textit{parallel 
in $M$} if there is an embedding of $S \times [0,1]$ in $M$ (called a 
\textit{parallelism from $S$ to $S\p$}) 
such that $S \times \{0\}=S$, $S \times \{1\}=S\p$, and $(\bd S) \times 
[0,1]$ lies in $\bd M$. If $S\p$ lies in $\bd M$ then $S$ is 
\textit{$\bd$-parallel in $M$}. The topological interior 
of $N$ in $M$ is denoted by $\Inte N$. 

\begin{lem} Let $W$ be a connected, non-compact 3-manifold which 
can be expressed as the union  
$W=\cup^{\infty}_{n=-\infty} X_n$ of compact, 
connected, \irr, \birr\ 3-manifolds $X_n$ such that 
$X_m \cap X_n=\emptyset$ 
for $|m-n|>1$ and $X_n \cap X_{n+1}=\bd X_n \cap \bd X_{n+1}$ is a 
compact, connected surface which is incompressible in $X_n$ and in 
$X_{n+1}$ and is not a disk. Then $\pi_1(W)$ is non-trivial and 
indecomposable. \end{lem}

\begin{proof} Standard arguments show 
that $Y_k=\cup^{k}_{n=-k} X_n$ is \irr\ and \birr. It follows that 
$\pi_1(Y_k)$ is non-trivial, non-infinite-cyclic, and indecomposable 
\cite[Theorem 5.2, Lemma 6.6]{He}. The incompressibility of each 
$X_n \cap X_{n+1}$ shows that 
$\pi_1(Y_k)$ injects into $\pi_1(M)$. We now apply Lemma 2.1. $\square$ 
\end{proof}

\section*{\normalsize 3. The Fox-Artin Arc}
\stepcounter{section}

\begin{thm} $\pi_1(S^3-A)$ is indecomposable, where $A$ is the 
Fox-Artin arc in Figure 1. \end{thm}

\begin{proof} Let $N$ be a tapered regular neighborhood of $A$. 
Thus $N$ is a 3-ball containing $A$ such that $A \cap \bd N = \bd A$, 
$A$ is isotopic in $N$ rel $\bd A$ to a diameter of $N$, and $N$ is 
tamely embedded in $S^3$ except at $\bd A$. Let $W=S^3-(\Inte N \cup 
\bd A)$. (We call $W$ the \textit{exterior} of $A$. We also use this 
term for the closure of the complement of a regular neighborhood of 
a tame submanifold of a manifold.)
Then $\pi_1(W) \cong \pi_1(S^3-A)$, and $\bd W=\bd N-\bd A$ is 
homeomorphic to an open annulus $S^1 \times \R$. It suffices to show 
that $W$ satisfies the hypotheses of Lemma 2.2. In the figures which 
follow we do not explicitly draw $N$, but its presence should be 
understood. 

$S^3-\bd A$ can be parametrized by $S^2 \times \R$ in such a way that 
$A$ meets each $S^2 \times [m, m+1]$, $m \in \Z$, in three arcs as 
indicated in Figure 3. 

\begin{figure}[h]
\begin{center}
\epsfig{file=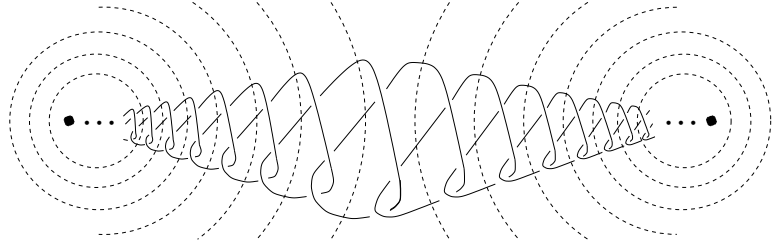, width=4in}
\end{center}
\caption{$S^3-\bd A$ parametrized as $S^2 \times \R$}
\end{figure}

It is natural to consider the exterior of the union of these three 
arcs in $S^2 
\times [m,m+1]$ and to regard $W$ as the union of these exteriors. 
Unfortunately these manifolds are cubes with two handles and so are 
not \birr. Instead we take $S^2 \times [2n-1,2n+1]$, $n \in \Z$, 
which also meets $A$ in three arcs, and let $X_n$ be the exterior 
of their union. 
The generic copy $X$ of $X_n$ is then the exterior of the union of 
the three arcs $\alpha$, $\beta$, and $\gamma$ in $S^2 \times [-1,1]$ 
as indicated in Figure 4. 

\begin{figure}[h]
\begin{center}
\epsfig{file=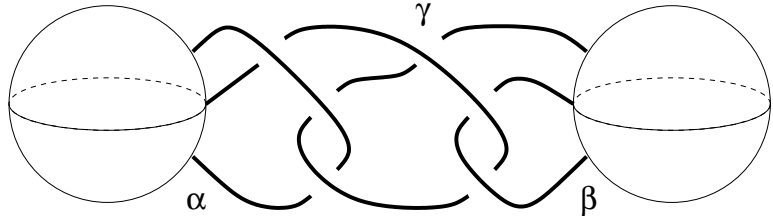, width=4in}
\end{center}
\caption{The arcs $\alpha$, $\beta$, and $\gamma$ in $S^2 \times 
[-1,1]$}
\end{figure}

Since no component of $X \cap (S^2 \times \{-1,1\})$  or of the closure 
of $\bd X-(S^2 \times \{-1,1\})$ is 
a disk it suffices to prove the following. 

\begin{lem} $X$ is \irr\ and \birr. \end{lem}

\begin{proof} Irreducibility follows from the Sch\"onflies theorem 
together with the fact that $X$ is a compact, connected submanifold 
of $S^3$ with connected boundary. 

The strategy for proving $\bd$-irreducibility is to exhibit $X$ as 
a double covering space of a solid torus $V$ branched over a certain 
properly embedded arc $\delta$ in $V$. If $\bd X$ were compressible, 
then by the $\Z_2$ case \cite{GL} of the equivariant loop theorem 
\cite{MY} there would be a compressing disk $\widetilde{D}$ for $\bd X$ 
such that either $\tau(\widetilde{D}) \cap \widetilde{D} = \nul$ or 
$\tau(\widetilde{D})=\widetilde{D}$, where $\tau$ is the non-trivial 
covering translation. Let $D$ be the image of $\widetilde{D}$ in $V$. 
In the first case $D$ would miss $\delta$. In the second case we could 
assume that $D$ would meet $\delta$ in a single transverse intersection 
point, since otherwise $\widetilde{D}$ would contain the fixed point 
set $\widetilde{\delta}$ of $\tau$, and we could reduce to the first case by 
replacing $\widetilde{D}$ by a nearby parallel disk. In both cases 
$D$ would be a compressing disk for $\bd V$ in $V$ since if $\bd D = 
\bd E$ for some disk $E$ in $\bd V$, then the preimage of $E$ in $X$ 
would have a component $\widetilde{E}$ with $\bd \widetilde{E} = 
\bd \widetilde{D}$. The proof is completed by showing that no such disk 
$D$ exists. 

By sliding one endpoint of each of $\alpha$ and of $\beta$ onto $\gamma$ 
we see that $X$ 
is homeomorphic to the exterior of the graph $\omega$ in $S^2 \times 
[-1,1]$  shown in Figure 5. 

\begin{figure}[h]
\begin{center}
\epsfig{file=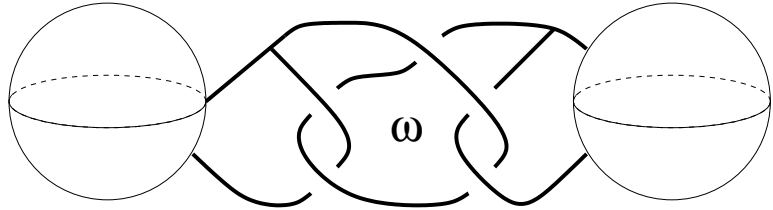, width=4in}
\end{center}
\caption{The graph $\omega$ in $S^2 \times [-1,1]$}
\end{figure}

This in turn is homeomorphic to the exterior $\widetilde{V}$ of the graph 
$\widetilde{\theta}$ in $S^3$ shown in Figure 6. 

\begin{figure}[h]
\begin{center}
\epsfig{file=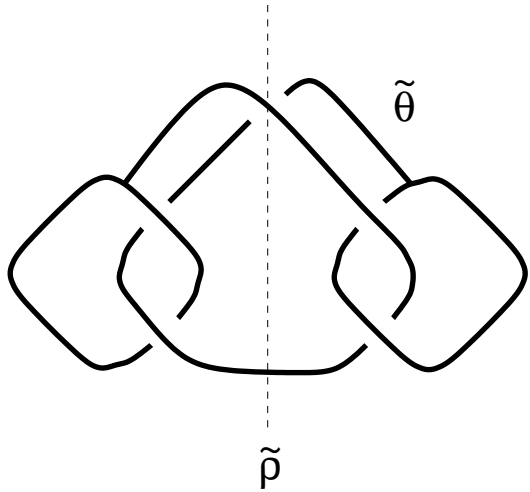, height=2.5in}
\end{center}
\caption{The graph $\widetilde{\theta}$ in $S^3$}
\end{figure}

This graph is invariant under the order two rotation $\tau$ about the 
simple closed curve $\widetilde{\rho}$. This involution defines 
a branched double covering $q:S^3 \rightarrow S^3$. The images $\theta$ 
and $\rho$ of $\widetilde{\theta}$ and $\widetilde{\rho}$ are 
shown in Figure 7.

\begin{figure}[h]
\begin{center}
\epsfig{file=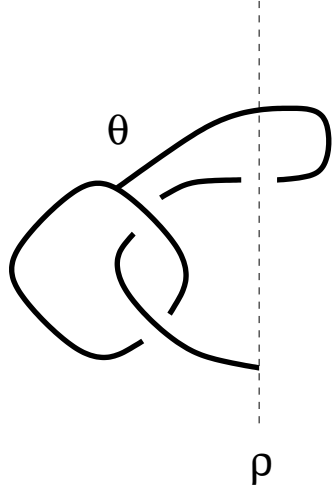, height=2.5in}
\end{center}
\caption{The quotient graph $\theta$ in $S^3$}
\end{figure}

Figure 8 shows a regular neighborhood $R$ of $\theta$ in $S^3$ and 
the arc $\delta=\rho \cap (S^3 - \Inte R)$. 
Figure 9 shows $R$ straightened by an isotopy to a standard solid 
torus.
Figure 10 moves the point at $\infty$ to a finite point. 
Figure 11 displays the solid torus $V=S^3-\Inte R$ containing $\delta$. 

\begin{figure}[h]
\begin{center}
\epsfig{file=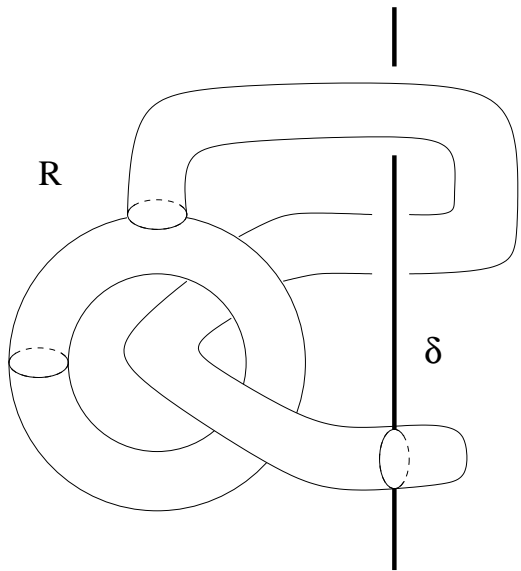, height=3in}
\end{center}
\caption{The regular neighborhood $R$ of $\theta$}
\end{figure}

\begin{figure}[h]
\begin{center}
\epsfig{file=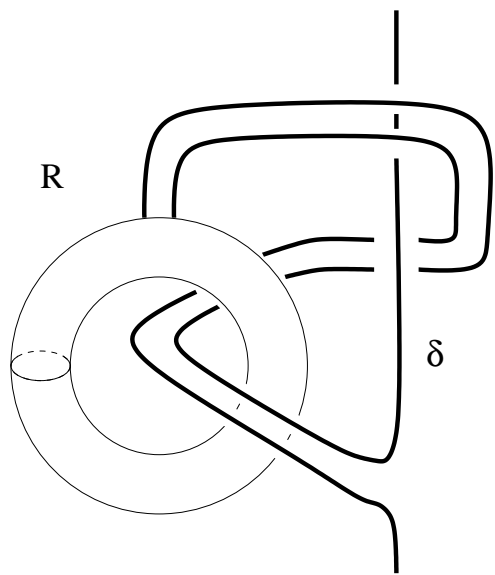, height=3in}
\end{center}
\caption{$R$ isotoped to a standard solid torus}
\end{figure}

\begin{figure}[h]
\begin{center}
\epsfig{file=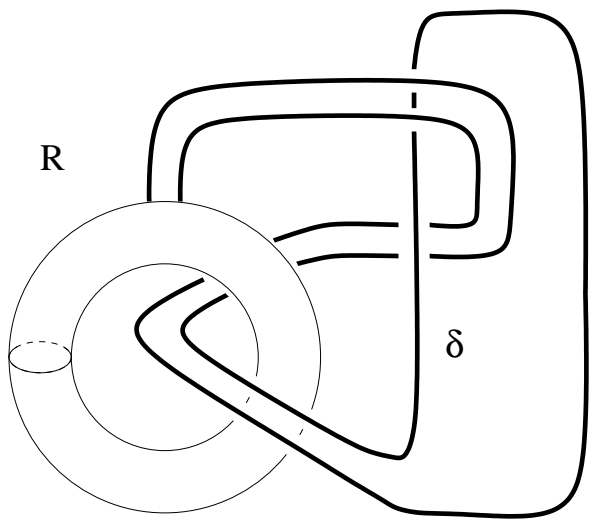, height=2.7in}
\end{center}
\caption{$\delta$ isotoped off the point at $\infty$}
\end{figure}

\clearpage

\begin{figure}[h]
\begin{center}
\epsfig{file=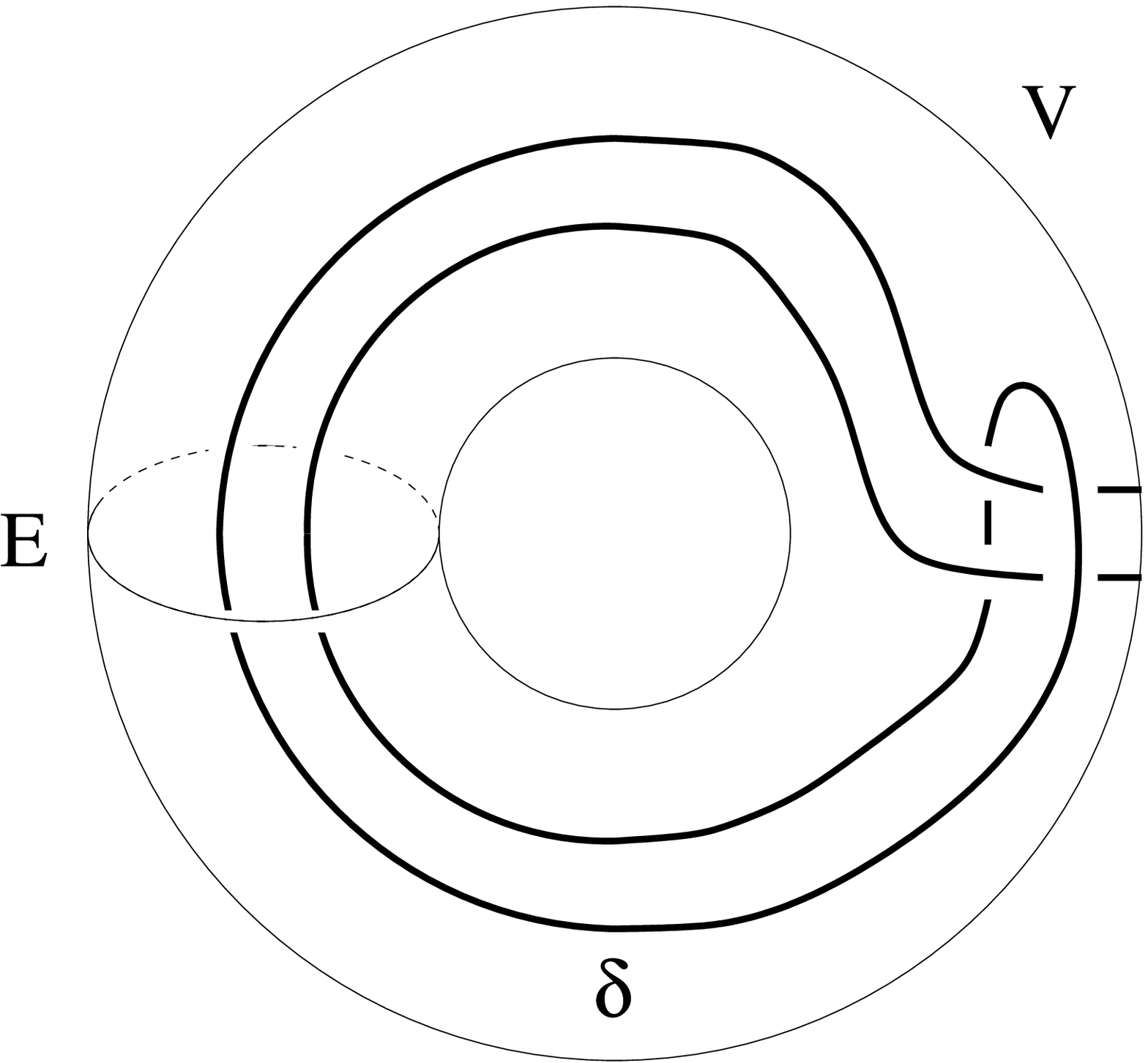, width=3in}
\end{center}
\caption{$\delta$ in $V=S^3-\Inte R$}
\end{figure}

\begin{lem} There is no meridinal disk $D$ in $V$ such that 
$D \cap \delta$ is either empty or a single transverse intersection 
point. \end{lem}

\begin{proof} Let $U$ be a regular neighborhood of $\delta$ in $V$. 
Let $F=\bd V - \Inte(U \cap \bd V)$ and $M=V - \Inte U$. It suffices to 
show that $F$ is incompressible in $M$ and that there is no properly 
embedded incompressible annulus $G$ in $M$ with one boundary component 
in the frontier (topological boundary) $C=Fr \,U$ of $U$ in $V$ and the 
other a curve in $F$ 
which bounds a meridinal disk $D$ in $V$ with $D \cap M = G$. 
Let $E$ be the meridinal disk shown in Figure 11. It meets $U$ 
in a pair of disks and so meets $M$ in a disk with two holes $S$. Let 
$V_0$ be the 3-ball obtained by splitting $V$ along $E$ and $M_0$ the 
3-manifold obtained by splitting $M$ along $S$. Then $E$ splits $\delta$ 
into three arcs $\delta_0$, $\delta_1$, and $\delta_2$, $U$ into the 
regular neighborhoods $U_0$, $U_1$, and $U_2$ of these arcs, $C$ into 
the three annuli $C_0$, $C_1$, and $C_2$, and $F$ into the surface 
$F_0$. See Figure 12. Let $S_0$ and $S_1$ be the copies of $S$ in $M_0$ 
which are identified to obtain $S$, where $S_0$ meets $C_0$ and $S_1$ 
meets $C_1$ and $C_2$. 

\begin{figure}[h]
\begin{center}
\epsfig{file=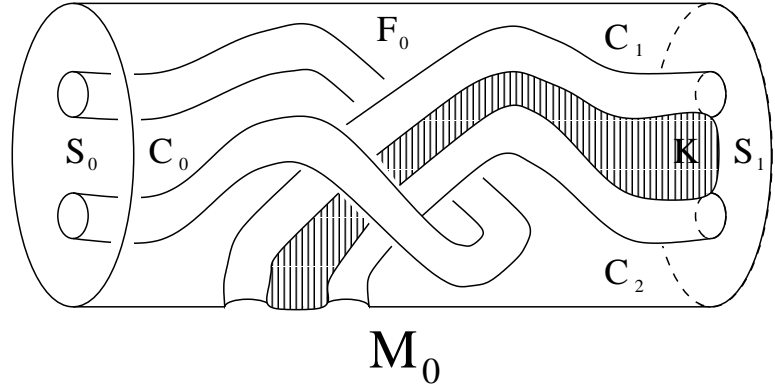, width=4in}
\end{center}
\caption{$M$ split along $S$ to obtain $M_0$}
\vspace{0.2in}
\end{figure}

Let $K$ be the disk in $M_0$ shown in Figure 12. Its boundary consists 
of one arc each in $F_0$, $S_1$, $C_1$, and $C_2$. Splitting $M_0$ along 
$K$ gives a 3-manifold $M_1$ which is homeomorphic to $(S_0 \cup C_0) 
\times [0,1]$ with $S_0 \cup C_0=(S_0 \cup C_0) \times \{0\}$. See 
Figure 13. $M_0$ is then obtained by attaching a 1-handle with cocore 
$K$ to $(S_0 \cup C_0) \times \{1\}$, so it is \irr. 

\begin{figure}[!h]
\begin{center}
\epsfig{file=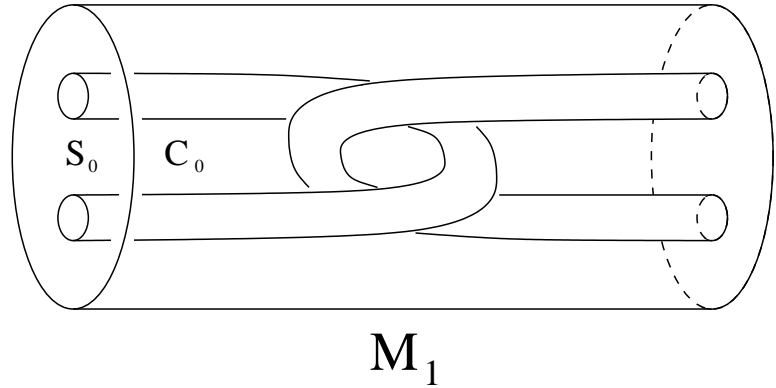, width=4in}
\end{center}
\caption{$M_0$ split along $K$ to obtain $M_1 \approx (S_0 \cup C_0) 
\times [0,1]$}
\end{figure}

We first show that $S$ is incompressible in $M$. 
It suffices to show that $S_0$ and $S_1$ are each 
incompressible in $M_0$. The first of these follows from our description 
above of $M_0$ as a product $I$-bundle with a 1-handle attached. 
The second follows from homology considerations. 

We next show that $F_0$ is incompressible in $M_0$.  
Suppose $L$ is a compressing disk. Then $\bd L$ separates 
one non-empty set of components of $\bd F_0$ from another. The seven 
possible partitions are all ruled out by a combination of  
homology arguments and the incompressibility of $S_0$. 

We now show that $S$ is \textit{$\bd$-incompressible rel $F$} in $M$. 
This means that 
whenever $L$ is a disk in $M$ such that $L \cap S$ is a properly 
embedded arc $\lambda$ in $S$ and $L \cap \bd M$ is an arc $\mu$ in 
$F$ such that $\lambda \cap \mu = \bd \lambda = \bd \mu$ and 
$\bd L=\lambda \cup \mu$, then there is an arc $\nu$ in 
$\bd S$ and a disk $L\p$ in $F$ such that $\mu \cap \nu = \bd \mu = 
\bd \nu$ and $\bd L\p=\mu \cup \nu$. 
It suffices to prove that $S_0$ and $S_1$ are $\bd$-incompressible 
rel $F_0$ in $M_0$. 

For $S_0$ this follows from homology considerations and the 
incompressibility of $F_0$ in $M_0$. For $S_1$ similar arguments 
reduce the problem to the case in which $\bd L=\lambda \cup \mu$ 
where $\lambda$ is an arc in $S_1$ such that $\bd \lambda$ lies in 
$S_1 \cap F_0$ and $\lambda$ separates $S_1 \cap C_1$ from $S_1 \cap 
C_2$ on $S_1$ and $\mu$ is an arc in $F_0$ separating $F_0 \cap C_1$ 
from $F_0 \cap C_2$. 

Isotop $L$ so that $K$ and $L$ are in general position and the arcs 
$K \cap S_1$ and $L \cap S_1$ meet in a single transverse intersection 
point. Then there is an arc $\xi$ in $K \cap L$ joining this point to a 
point in $K \cap F_0$. Since $M_0$ is \irr\ we may assume that in 
addition $K \cap L$ contains no simple closed curves. The intersection 
then consists of $\xi$ and possibly some arcs $\eta$ with $\bd \eta$ in 
$K \cap F_0$. Assume $\eta$ is outermost on $L$. Let $\zeta$ be an arc 
in $\bd L$ such that $\zeta \cup \eta$ bounds a disk $L_0$ in $L$ whose 
interior misses $K$. Let $\varepsilon$ be the arc on $K \cap F_0$ with 
$\bd \varepsilon = \bd \eta = \bd \zeta$. There is a disk $K_0$ in $K$ 
such that $\bd K_0=\eta \cup \varepsilon$. Then $K_0 \cap L_0=\eta$ and 
$K_0 \cup L_0$ is a disk with boundary $\zeta \cup \varepsilon$. Since 
$F_0$ is incompressible in $M_0$ this curve bounds a disk $F_1$ in 
$F_0$. Since $M_0$ is \irr\ $K_0 \cup L_0 \cup F_1$ bounds a 3-ball 
$B_0$ in $M_0$. Note that $\xi \cap B_0=\nul$. An isotopy of $L$ which 
moves $L_0$ across $B_0$ to $K_0$ and then off $K_0$ removes $\eta$ and 
possibly other components of $K \cap L$ but does not affect $\xi$. 

Thus we may assume that $K \cap L=\xi$. We now split $M_0$ along $K$ to 
obtain $M_1$, as before. This splits $L$ into disks $L_0$ and $L_1$ 
either of which we can take as a compressing disk for $(S_0 \cup C_0) 
\times \{1\}$ in $M_1=(S_0 \cup C_0) \times [0,1]$. This contradiction 
completes the proof that $S_0$ and $S_1$ are $\bd$-incompressible rel 
$F_0$ in $M_0$ 
and hence that $S$ is $\bd$-incompressible rel $F$ in $M$. 

Now suppose that $D$ is a compressing disk for $F$ in $M$. Put $D$ in 
general position with respect to $S$ so that $D \cap S$ has a minimal 
number of components. By the incompressibility of $S$ and the 
irreducibility of $M$ none of them are simple closed curves. Since 
$S$ is $\bd$-incompressible rel $F$ in $M$ none of them can be arcs, so 
$D \cap S=\nul$. Since $F_0$ is incompressible in $M_0$ we have that 
$D$ cannot exist. 

Finally suppose that $G$ is an incompressible annulus in $M$ with one 
boundary component in $C$ and the other a curve in $F$ which bounds 
a meridinal disk $D$ of $V$ such that $D \cap M = G$. 
We may assume that the first  
boundary component misses $S$, that $G$ is in general position with 
respect to $S$ and that among all such annuli in its isotopy class 
$G \cap S$ has a minimal number of components. Then none of these 
components is a simple closed curve which bounds a disk in $S$ or 
in $G$ or is an arc joining the two components of $\bd G$. 

Suppose some component $\kappa$ of $G \cap S$ is a simple closed curve. 
Then we may assume that $\kappa$ and $G \cap C$ form the boundary of 
a subannulus $G_0$ of $G$ which lies in $M_0$. If $\kappa$ lies in $S_0$, 
then for homological reasons $G \cap C$ must lie in $C_0$. We can isotop 
$G_0$ so that it misses $K$. Hence $G_0$ lies in $M_1=(S_0 \cup C_0) 
\times [0,1]$. By \cite[Corollary 3.2]{Wa} $G_0$ is parallel to an annulus 
in $(S_0 \cup C_0) \times \{0\}$ and so $\kappa$ can be removed by an 
isotopy, contradicting minimality. If $\kappa$ lies in $S_1$, then for 
homological reasons $G \cap C$ must be in $C_1$ or $C_2$, say $C_1$. 
Let $M_2=M_0 \cup U_0$. Then $M_2$ is homeomorphic to $S_1 \times 
[0,1]$ with $S_1=S_1 \times \{1\}$. Now $G_0$ is incompressible in 
$M_2$ and can be isotoped keeping $\kappa$ fixed to an annulus $G\p_0$ 
such that $\bd G\p_0$ lies in $S_1$. It then follows from  
\cite[Corollary 3.2]{Wa} that $G\p_0$ is parallel to an annulus in $S_1$ and hence 
$G_0$ is $\bd$-parallel in $M_2$. Since this parallelism does not 
meet $U_0$ we have that $G_0$ is $\bd$-parallel in $M_0$. It follows 
that $\kappa$ can be removed by an isotopy, again contradicting 
minimality. 

Hence any component of $G \cap S$ must be an arc whose boundary lies 
in $F \cap S$. Since $S$ is $\bd$-incompressible rel $F$ in $M$ and 
$S$ is incompressible in $M$ any outermost such arc can be removed by 
an isotopy. Thus $G \cap S=\nul$, and we may regard $G$ as lying in 
$M_0$. For homological reasons $G \cap C$ must lie in $C_1$ or $C_2$, 
say $C_1$. Since $D$ is a 
meridinal disk of $V$ we must have for homological reasons that $\bd D$ 
splits $F_0$ into two components such that one contains $F_0 \cap S_0$ 
and $F_0 \cap C_1$ and the other contains $F_0 \cap S_1$ and 
$F_0 \cap C_2$. Let 
$M\p_1=M_0 \cup U_1$. Then $M\p_1$ is homeomorphic to $(S_0 \cup C_0) 
\times [0,1]$ with $S_0 \cup C_0=(S_0 \cup C_0) \times \{0\}$. 
So $D$ is a compressing disk for $\bd M\p_1 - (S_0 \cup C_0)$ in 
$M\p_1$. This contradiction completes the proof of Lemma 3.3. $\square$

\end{proof}

This completes the proof of Lemma 3.2. $\square$ \end{proof}

This completes the proof of Theorem 3.1. $\square$ \end{proof}

\section*{\normalsize 4. Embeddings of Torus Knot Groups}
\stepcounter{section}

In this section we prove a technical result concerning embeddings 
of torus knot groups in torus knot groups which will be used in 
the next section to distinguish among the fundamental groups of 
the complements of a certain uncountable collection of arcs. 
Recall that the fundamental group of the complement of a $(p,q)$ 
torus knot is the group $G_{p,q}=\langle x,y\,|\,x^p=y^q \rangle$.  

\begin{lem} Let $p$, $q$, $r$, and $s$ be primes such that 
$p<q$ and $r<s$. Then $G_{p,q}$ embeds in $G_{r,s}$ if and only 
if $p=r$ and $q=s$. \end{lem}

\begin{proof} Let $Z(G)$ denote the center of the group $G$. 
Recall that $Z(G_{p,q})$ is an infinite cyclic group generated by 
$x^p$ and that $G/Z(G_{p,q})\cong \Z_p*\Z_q$. Recall also that a 
free product of two non-trivial groups has trivial center and that 
any element of finite order in a free product is conjugate to an 
element of one of the factors. (See \cite[pp. 140--141, 100--101]{Ma}.) 

We may assume that $G_{p,q}$ is a subgroup of $G_{r,s}$. Let 
$K=G_{p,q} \cap Z(G_{r,s})$. Then $K$ is a subgroup of $Z(G_{p,q})$ 
and is the kernel of the restriction of the natural projection 
$G_{r,s} \rightarrow \Z_r*\Z_s$ to $G_{p,q}$. If $u \in G_{r,s}$, 
then let $\bar{u}$ denote its image in $\Z_r*\Z_s$. 

Suppose $K=Z(G_{p,q})$. Then we have an embedding 
$\Z_p*\Z_q \rightarrow \Z_r*\Z_s$. Since $\bar{x}$ has order $p$ it 
must be conjugate to an element of $\Z_r$ or of $\Z_s$, hence 
$p|r$ or $p|s$, hence since $r$ and $s$ are prime we have $p=r$ or 
$p=s$. Similarly $q=r$ or $q=s$. Since $p<q$ and $r<s$ we must 
have $p=r$ and $q=s$. 

Now suppose that $K$ is a proper subgroup of $Z(G_{p,q})$. Then it 
is generated by $x^{pk}$ for some $k \geq 0$, $k \neq 1$. Let 
$G_{p,q,k}=G_{p,q}/K$. It embeds in $\Z_r*\Z_s$ and has presentation 
$\langle\bar{x}, \bar{y} \,|\, \bar{x}^p=\bar{y}^q, \bar{x}^{pk}=1\rangle$. 
By the Kurosh subgroup theorem \cite{Ku,Ma} $G_{p,q,k}$ must be 
a free product of cyclic groups and so must either be cyclic or 
have trivial center. It thus suffices to show that neither of these 
is the case. 

For $k=0$ this group is just $G_{p,q}$, and we are done. So assume 
$k \geq 2$. Define functions $f,g:\Z_{pqk}\rightarrow\Z_{pqk}$ by 
$f(n)=n+q \bmod pqk$ and $g(n)=n+p \bmod pqk$. Then 
$f$ and $g$ are one to one and so may be regarded as elements of the 
symmetric group $\mathcal{S}_{pqk}$. Define $\psi:G_{p,q,k} \rightarrow 
\mathcal{S}_{pqk}$ 
by $\psi(\bar{x})=f$ and $\psi(\bar{y})=g$. Then $\psi$ is well defined 
because $f^p(n)=n+pq=n+qp=g^q(n)$ and $f^{pk}(n)=n+pkq = n 
\bmod pqk$. Since $\psi(\bar{x}^p)=f^p \neq id$ we have that 
$Z(G_{p,q,k})$ is non-trivial. Since $G_{p,q,k}$ maps onto $\Z_p*\Z_q$ 
it is non-cyclic, and so we are done. $\square$ \end{proof}

\section*{\normalsize 5. Uncountably Many Arcs}
\stepcounter{section}

\begin{thm} There are uncountably many arcs $A_i$ in $S^3$ such 
that:
\begin{enumerate}
\item $\pi_1(S^3-A_i)$ is indecomposable and non-trivial.
\item $\pi_1(S^3-A_i)$ and $\pi_1(S^3-A_j)$ are isomorphic if and 
only if $i=j$. 
\item $A_i$ is wildly embedded precisely at its endpoints. 
\end{enumerate}
\end{thm}

\begin{proof} We first outline the proof and then fill in the details 
with a sequence of lemmas. 

The construction of the $A_i$ will have a pattern similar to that of  
the Fox-Artin arc. $S^3-\bd A_i$ will be parametrized as $S^2 \times \R$, 
and for each integer $n$ we will have that $A_i$ meets $S^2 \times 
[n,n+1]$ in three properly embedded arcs $\alpha_n$, $\beta_n$, and 
$\gamma_n$, where $\alpha_n$ runs from $S^2 \times \{n\}$ to itself, 
$\beta_n$ runs from $S^2 \times \{n+1\}$ to itself, and $\gamma_n$ runs 
from $S^2 \times \{n\}$ to $S^2 \times \{n+1\}$. These arcs will be 
chosen so that the exterior $X_n$ of $\alpha_n \cup \beta_n \cup 
\gamma_n$ in $S^2 \times [n,n+1]$ is \irr\ and \birr. Hence by  
Lemma 2.2 we will have that $\pi_1(S^3-A_i)$ is indecomposable 
and non-trivial. Thus $A_i$ will be wild. It will clearly 
be tame at points not in $\bd A_i$. It will be wild at both endpoints 
since otherwise its complement would be simply connected. (Any meridian  
of the arc would bound a disk consisting of an annulus which follows 
the arc to a tame endpoint and is then is capped off by a disk behind it. 
In fact it can be shown as in \cite[Example 1.2]{FA} that $S^3-A_i$ would be 
homeomorphic to \RRR.) 

A map is $\pi_1$-\textit{injective} if it induces an injection on 
fundamental groups; the same term is applied to a submanifold if 
its inclusion map has this property. 
The arcs will be chosen so that the interior of $X_n$ will contain a 
$\pi_1$-injective submanifold $Q_n$ which is homeomorphic to the exterior 
of a $(p_n,q_n)$ torus knot in $S^3$, where $p_n$ and $q_n$ are primes 
with $p_n < q_n$. It will follow 
from the $\bd$-irreducibility of all the $X_m$ that $\pi_1(S^3-A_i)$ 
will have a subgroup isomorphic to $\pi_1(Q_n)$. Moreover it will be 
shown that any subgroup of $\pi_1(S^3-A_i)$ which is isomorphic to a 
$(p,q)$ torus knot group for primes $p$ and $q$ with $p < q$ must be 
isomorphic to one of the $\pi_1(Q_n)$. We then let $J$ be the set of 
all pairs of primes $(p,q)$ with $p < q$ and let $2^J$ be the set of 
all subsets of $J$. For each non-empty $i \in 2^J$ we 
construct an arc $A_i$ as 
above such that the $(p,q)$ torus knot subgroups of $\pi_1(S^3-A_i)$ 
with $(p,q) \in J$ are precisely those for which $(p,q) \in i$. It 
follows that $\pi_1(S^3-A_i)$ and $\pi_1(S^3-A_j)$ are isomorphic 
if and only if $i=j$. Since $2^J$ is uncountable we will be done. 

We next recall some terminology. Let $M$ be a compact, connected,  
orientable 3-manifold. We say that $M$ is \textit{atoroidal} 
if every properly embedded, incompressible torus $S^1 \times S^1$ in $M$ 
is $\bd$-parallel in $M$ and is \textit{anannular} if every 
properly embedded, incompressible annulus $S^1 \times [0,1]$ in $M$ is 
$\bd$-parallel in $M$. If $M$ is \irr, \birr, anannular and atoroidal, 
contains a 2-sided, properly embedded incompressible surface, and is not 
a 3-ball, then $M$ is \textit{excellent}; the 
same term is applied to a compact, properly embedded 1-manifold in 
a compact 3-manifold $P$ if its exterior in $P$ has these properties.

\begin{lem} Let $Y\p$ and $Y\pp$ be excellent 3-manifolds. Suppose 
$Y=Y\p \cup Y\pp$, where $S=Y\p \cap Y\pp 
=\bd Y\p \cap \bd Y\pp$ is a compact surface such that $S$ is incompressible 
in $Y\p$ and in $Y\pp$, $\bd Y\p-\Inte S$ is incompressible in $Y\p$, 
$\bd Y\pp-\Inte S$ is incompressible in $Y\pp$, and each component of 
$S$ has negative Euler characteristic. Then $Y$ is excellent. \end{lem} 

\begin{proof} This is \cite[Lemma 2.1]{My excel}. $\square$ \end{proof} 

We now construct the arcs. Let $R$ be an unknotted solid torus in the 
interior of $S^2 \times [0,1]$. 
Let $P=S^2 \times [0,1]-\Inte R$. (We say that $R$ 
is \textit{unknotted} if there is a properly embedded disk $E$ in $P$ 
such that $\bd E \subseteq \bd R$ and a meridinal disk $D$ of $R$ such 
that $\bd D$ and $\bd E$ meet transversely in a single point.)  

\begin{lem} There exist disjoint properly embedded arcs $\alpha$, $\beta$, 
and $\gamma$ in $P$ such that $\bd \alpha \subseteq S^2 \times \{0\}$, 
$\bd \beta \subseteq S^2 \times \{1\}$, $\gamma$ has one endpoint 
in $S^2 \times \{0\}$ and the other in $S^2 \times \{1\}$, and 
$\alpha \cup \beta \cup \gamma$ is excellent. \end{lem}

\begin{proof} Let $\alpha\p$, $\beta\p$, and $\gamma\p$ be any arcs 
in $P$ whose boundaries satisfy the given conditions. By 
\cite[Theorem 1.1]{My excel} any compact, properly embedded 1-manifold in a compact, 
connected, orientable 3-manifold which meets each 2-sphere boundary 
component in at least two points is homotopic relative its boundary 
to a properly embedded 1-manifold which is excellent. Let $\alpha$, 
$\beta$, and $\gamma$ be the 
respective components of this new 1-manifold. 

For those who prefer a more concrete construction of such arcs we 
give an alternative proof at the end of this section. $\square$ 
\end{proof} 

Now let $Q$ be the exterior of a $(p,q)$ torus knot in $S^3$, where 
$(p,q) \in J$. Glue $P$ and $Q$ together by identifying $\bd R$ with 
$\bd Q$ in such a way that $\bd E$ is identified with a meridian of 
$\bd Q$. Then the union of $Q$ and a regular neighborhood of $E$ in 
$P$ is a 3-ball, and so $P \cup Q$ is homeomorphic to $S^2 \times 
[0,1]$. Let $Y$ be the exterior of $\alpha \cup \beta \cup \gamma$ in 
$P$ and $X=Y \cup Q$. It follows from the irreducibility 
and $\bd$-irreducibility of $Y$ and of $Q$ that $X$ is \irr\ and \birr\  
and that $Q$ is $\pi_1$-injective in $X$. 

We now repeat this construction using $(p_n,q_n)$ torus knots with 
$(p_n,q_n)\in i$ to obtain $\alpha_n$, $\beta_n$, 
$\gamma_n$, $P_n$, $Q_n$, $Y_n$, and $X_n$ contained in $S^2 \times 
[n,n+1]$. We construct an arc $A_i$ by identifying the endpoints of 
the arcs so that the arcs occur in the sequence $\ldots, \gamma_n, 
\alpha_{n+1}, \beta_n, \gamma_{n+1}, \ldots$ on $A_i$. The exterior $W_i$ 
of $A_i$ then satisfies the hypotheses of Lemma 2.2, and so 
$\pi_1(S^3-A_i)$ is indecomposable and non-trivial. Moreover the 
incompressibility of each $X_n \cap X_{n+1}$ implies that each $Q_n$ 
is $\pi_1$-injective in $W_i$. 

We next review some characteristic submanifold theory \cite{Ja,JS,Jo}, 
following \cite{JS} but restricting attention to the special case 
which we will need. We first refine our notion of parallel 
surfaces. 
A pair $(M,F)$ is an \textit{\irr\ 3-manifold pair} if $M$ is a compact, 
orientable, \irr\ 3-manifold and $F$ is a compact, incompressible surface 
in $\bd M$. Let $S$ and $S\p$ be disjoint compact surfaces in $M$ such 
that $S$ is properly embedded in $M$, $S\p$ is either properly embedded 
in $M$ or contained in $\bd M$, and $\bd S \cup \bd S\p$ is contained 
in $F$. We say that $S$ and $S\p$ are \textit{parallel in $(M,F)$} 
if there is a parallelism $S \times [0,1]$ from $S$ to $S\p$ such that 
$(\bd S) \times [0,1]$ is contained in $F$; if $S\p \subseteq F$ we say 
that $S$ is \textit{$F$-parallel}. 
Our old definitions of 
``parallel'' and ``$\bd$-parallel'' in $M$ correspond to the case of 
$F=\bd M$.  

The \textit{characteristic pair} of the \irr\ 3-manifold pair $(M,\bd M)$ 
is a certain \irr\ 3-manifold pair $(\Sigma,\Phi)$ such that $\Sigma 
\subseteq M$ and $\Sigma \cap \bd M=\Phi$.  
For its definition and proof of existence see \cite[Chapter V]{JS}. 
We will limit our discussion to two basic issues: using $(\Sigma,\Phi)$ 
and recognizing $(\Sigma,\Phi)$. The property we will use is 
that any $\pi_1$-injective map from a 
Seifert fibered space with non-cyclic fundamental group into $M$ 
which is not homotopic to a map whose image lies in $\bd M$ must be 
homotopic to a map whose image lies in $\Sigma$ \cite[p. 138]{JS}. 

We will recognize $\Sigma$ by recognizing its components and using 
the Splitting Theorem \cite[p. 157]{JS} to recognize the frontier 
$Fr \,\Sigma$ of $\Sigma$ in $M$. The components $(\sigma,\varphi)$ 
of $(\Sigma,\Phi)$ are \textit{Seifert pairs}, i.e.\ $\sigma$ is either 
an $I$-bundle over a compact surface with $\varphi$ the associated 
$\bd I$-bundle or $\sigma$ is a Seifert fibered space with $\varphi$ a 
union of fibers in $\bd \sigma$. One of the properties we will need is 
that the inclusion map 
from $(\sigma,\varphi)$ into $(M,\bd M)$ is not homotopic as a map 
of pairs to a map whose image lies in $\Sigma-\sigma$. Also  
the components of $Fr \,\Sigma$ are incompressible annuli and tori 
none of which is $\bd$-parallel in $M$ though some components may be 
parallel in $(M,\bd M)$ to each other. (See the examples in  
\cite[Chapter IX]{Ja}.) A union $Fr^*\,\Sigma$ of components of 
$Fr \,\Sigma$ 
such that no two components of $Fr^*\, \Sigma$ are parallel in 
$(M,\bd M)$ to each other and $Fr^*\, \Sigma$ is maximal with respect 
to inclusion among all such unions is called a \textit{reduction} of 
$Fr \,\Sigma$. We call the components of $Fr \,\Sigma-Fr^* \,\Sigma$ 
\textit{redundant} components of $Fr \,\Sigma$. 
Now suppose we are given a compact, properly embedded 
surface $\mathcal{T}$ in $M$ satisfying the following two conditions: 
\renewcommand{\theenumi}{\alph{enumi}}
\renewcommand{\labelenumi}{(\theenumi)}
\begin{enumerate}
\item The components of $\mathcal{T}$ are incompressible 
annuli and tori none of which is $\bd$-parallel in $M$. 
\item Let $(M\p,\bd\bp\!M)$ be the pair 
obtained by splitting $M$ along $\mathcal{T}$ and $\bd M$ along 
$\bd \mathcal{T}$. Then each component $(N,L)$ of $(M\p,\bd\bp\! M)$ is 
either a Seifert pair or a \textit{simple pair}, i.e. every 
incompressible, properly embedded torus in $N$ or annulus in $N$ with 
boundary in $\Inte L$ is either $L$-parallel or parallel in $(N,L)$ 
to a component of $\bd N-\Inte L$. 
\end{enumerate}
If $\mathcal{T}$ is minimal with respect to inclusion among all 
compact, properly embedded surfaces in $M$ satisfying (a) and (b), 
then by the Splitting Theorem $\mathcal{T}$ 
is isotopic to $Fr^* \,\Sigma$. 

Now let $M_k=\cup_{n=-k}^k X_n$ and $C_k=\cup_{n=-k}^k Q_n$.

\begin{lem} $(M_k,\bd M_k)$ is an \irr\ 3-manifold pair, and 
its characteristic pair $(\Sigma,\Phi)=(C_k,\nul)$. \end{lem}

\begin{proof} The irreducibility and \bd-irreducibility of $M_k$ 
and the incompressibility of $\bd C_k$ in $M_k$ follow from the 
irreducibility and \bd-irreducibility of the $X_n$, the 
incompressibility of the $X_n \cap X_{n+1}$ in $X_n$ and in $X_{n+1}$, 
and the incompressibility of $\bd Q_n$ in $X_n$. 

Let $\mathcal{T}=\bd C_k$. Since $\bd M_k$ is a surface of genus two 
no component of $\mathcal{T}$ is $\bd$-parallel in $M_k$. 
The components of $(M\p_k,\bd\bp\! M_k)$ 
are the $(Q_n,\nul)$ and $(Z,\bd M_k)$, where $Z=\cup_{n=-k}^k Y_n$. 
Each $Q_n$ is a Seifert fibered space. By Lemma 5.2 we have that 
$Z$ is excellent 
and therefore $(Z,\bd M_k)$ is a simple pair. Thus $\mathcal{T}$ 
satisfies properties (a) and (b). Deleting any components of 
$\mathcal{T}$ gives a surface which splits $M_k$ into components 
one of which, say $N$, is the union of $Z$ and some of the $Q_n$. 
Now $N$ is not Seifert fibered since it contains $\bd M_k$. 
It is not an $I$-bundle over a compact surface $S$ since $S$ would be 
covered by $\bd M_k$, and so $\pi_1(S)\cong\pi_1(N)$  
could not contain the $\mathbf{Z}\oplus\mathbf{Z}$ subgroup 
$\pi_1(\bd Q_n)$. Finally $(N,\bd M_k)$ is not a simple pair because 
$\bd Q_n$ is not \bd-parallel in $N$. Thus $\mathcal{T}$ is minimal 
with respect to inclusion among surfaces satisying (a) and (b). 
So by the Splitting Theorem $\mathcal{T}=Fr^* \,\Sigma$. 

By arguments similar to those applied above to $N$ we have that 
$(Z,\bd M_k)$ is not a Seifert pair. So if there are no redundant 
components we must have $(\Sigma,\Phi)=(C_k,\nul)$, and we are done. 

Suppose there is a redundant component. Then it must be  a torus which 
is parallel in $(M_k,\bd M_k)$ to $\bd Q_n$ for some $n$; denote 
it by $T_n$. Thus we may assume that there is an embedding of 
$T_n \times [0,1]$ in $M_k$ such that $T_n \times [0,1]$ meets $Q_n$ in 
$T_n \times \{0\}=\bd Q_n$, $T_n \times \{1\}=T_n$, and 
$T_n \times (0,1)$ contains all other redundant tori which are parallel 
to $\bd Q_n$. 
If there are such extra redundant tori, then they are isotopic in 
$T_n \times [0,1]$ to tori of the form $T_n \times \{t\}$ 
\cite[Corollary 3.2]{Wa}. 
It follows that there is some component $\sigma$ of $\Sigma$ of the form 
$T_n \times [r,s]$. Its inclusion map into $M_k$ is homotopic to a map 
whose image lies in $\Sigma-\sigma$, contradicting one of the properties 
of $\Sigma$. 

Thus there are no extra redundant tori. Now let $Z\p$ be 
the closure of the complement in $Z$ of the union of all the products 
$T_n \times [0,1]$. Then $Z\p$ is homeomorphic to $Z$, and so 
$(Z\p,\bd M_k)$ is a simple pair which is not Seifert pair. Thus 
$T_n \times [0,1]$ is a component of $\Sigma$, and $(Q_n,\nul)$ is 
a simple pair. Now in fact $(Q_n,\nul)$ actually \textit{is} a simple 
pair. However, it is also a Seifert fibered space with non-cyclic 
fundamental group. Its inclusion map cannot be homotopic to a map 
whose image lies in $\bd M_k$ because $\pi_1(M_k)$ has no 
$\mathbf{Z}\oplus\mathbf{Z}$ subgroups. Thus it must be homotopic to 
a map whose image lies in some component $\sigma$ of $\Sigma$. 
In particular the image lies in the complement of $Q_n$. 

Now it follows 
from \cite[Squeezing Theorem, p. 139]{JS} or \cite[Theorem IX.12]{Ja} 
that $Q_n$ is actually isotopic to a submanifold of $\sigma$. This 
fact can be used to contradict our knowledge of the structure of $Z\p$. 
We choose, however, to give the following somewhat more direct argument. 

Let $p:\widetilde{M}_k \rightarrow 
M_k$ be the covering map corresponding to $\pi_1(Q_n)$. There is a 
component $\widetilde{Q}_n$ of $p^{-1}(Q_n)$ such that the restriction 
$\widetilde{Q}_n \rightarrow Q_n$ of $p$ is a homeomorphism and  
$\pi_1(\widetilde{Q}_n) \rightarrow \pi_1(\widetilde{M}_k)$ is an 
isomorphism. It follows that $\pi_1(\bd \widetilde{Q}_n) \rightarrow 
\pi_1(\widetilde{M}_k - \Inte \widetilde{Q}_n)$ is an isomorphism. 
Now the homotopy of $Q_n$ into its complement lifts to a homotopy of 
$\widetilde{Q}_n$ into $\widetilde{M}_k - \Inte \widetilde{Q}_n$.  
This implies that $\pi_1(Q_n)$ is abelian, which is not the case. 
$\square$ \end{proof} 

We now suppose that $\pi_1(S^3-A_i)$ and $\pi_1(S^3-A_j)$ are 
isomorphic. Then $\pi_1(W_i)$ and $\pi_1(W_j)$ are isomorphic, 
where $W_i$ and $W_j$ are the exteriors of $A_i$ and $A_j$, 
respectively. Since these spaces are \irr\ and orientable, the 
sphere theorem implies that they are aspherical. Hence there is 
a map $h:W_j \rightarrow W_i$ such that $h_*:\pi_1(W_j) \rightarrow 
\pi_1(W_i)$ is an isomorphism. We then restrict $h$ to a $(p,q)$   
torus knot space arising in the construction of $A_j$. This map is 
$\pi_1$-injective. Its image lies in some $M_k$. Since $\pi_1(\bd M_k)$ 
has no $\mathbf{Z}\oplus\mathbf{Z}$ subgroups Lemma 5.4 implies that 
it is homotopic to a map whose image lies in some $(r,s)$ torus 
knot space arising in the construction of $A_i$. By Lemma 4.1 we 
have that $(p,q)=(r,s)$. Thus $j \subseteq i$. The symmetric 
argument shows that $i \subseteq j$, concluding the proof of 
Theorem 5.1. $\square$ \end{proof}

\begin{figure}[h]
\begin{minipage}{3in}
\begin{center}
\epsfig{file=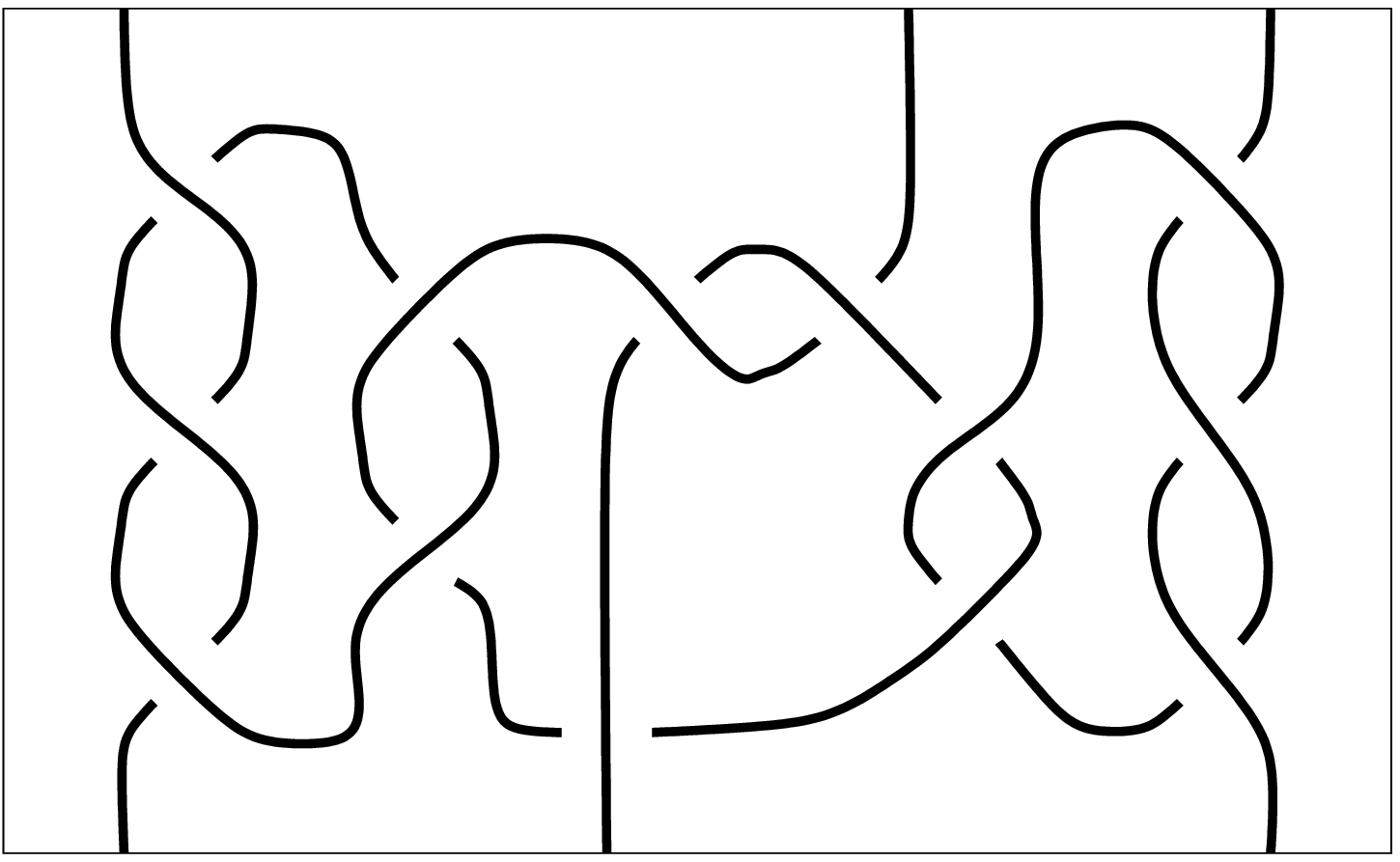, height=1.8in}
\end{center}
\caption{An excellent 3-tangle}
\end{minipage}
\begin{minipage}{2.1in}
\begin{center}
\epsfig{file=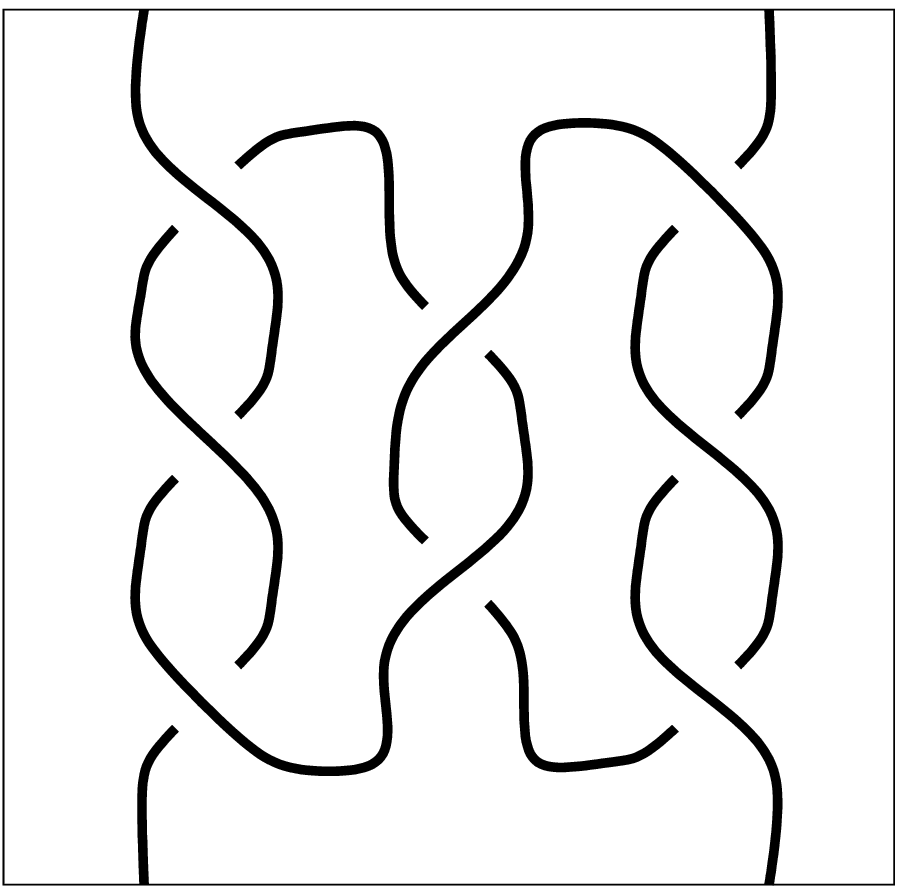, height=1.8in}
\end{center}
\caption{An excellent 2-tangle}
\end{minipage}
\end{figure}

\begin{proof}[\textnormal{\textbf{Alternative Proof of Lemma 5.3}}] 
Figure 14 shows a three component tangle in a 3-ball. Figure 15 
shows a two component tangle in a 3-ball. 
By \cite[Proposition 4.1]{My homology} and  
\cite[Proposition 4.1]{My simple} these two tangles are excellent. 
Let $Y\p$ and $Y\pp$ be their respective exteriors.  

We glue $Y\p$ and $Y\pp$ together as indicated in Figure 16 to obtain 
the exterior $Y$ of the union of the arcs $\alpha$, $\beta$, and $\gamma$ 
in the space $P$ obtained by removing the interior of an unknotted 
solid torus $R$ contained in the interior of $S^2 \times [0,1]$. 
$S=Y\p \cap Y\pp=\bd Y\p \cap \bd Y\pp$ has two components; 
each is a disk with two holes. Since a compact 
surface contained in an incompressible boundary component of a 
compact 3-manifold is incompressible if none of the components of 
its complement in the boundary component has closure a disk, we have 
that $S$ is incompressible in $Y\p$ and in $Y\pp$. We now apply Lemma 5.2 
to conclude that $Y$ is excellent.  
$\square$ \end{proof}

\begin{figure}[!h]
\begin{center}
\epsfig{file=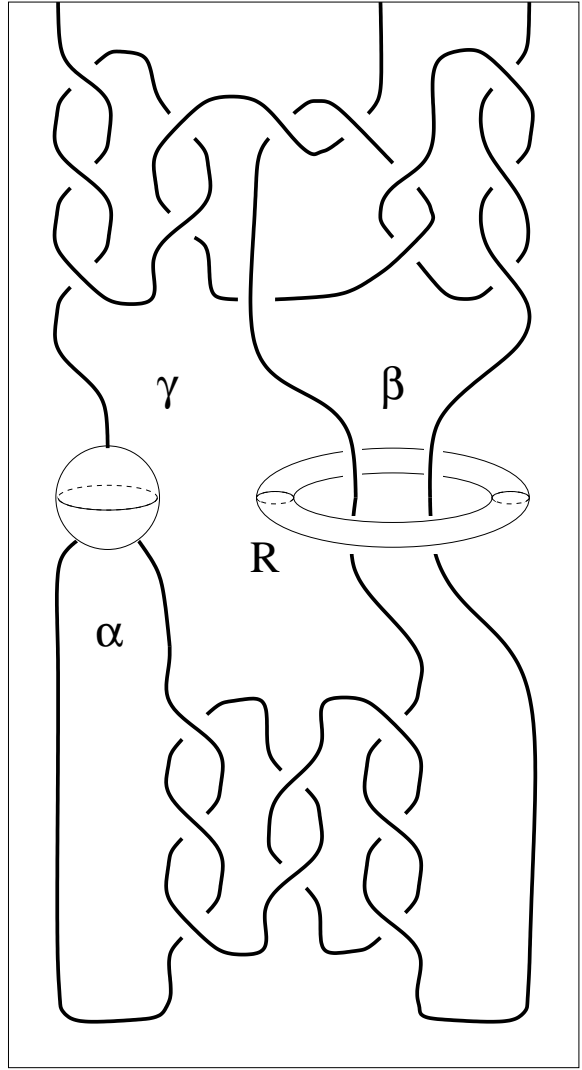, width=3in}
\end{center}
\caption{The three arcs $\alpha$, $\beta$, and $\gamma$ in 
$(S^2 \times [0,1])-\Inte R$}
\end{figure}


\end{document}